\documentclass[11pt]{amsart}
\usepackage{ifthen, citesort}
\usepackage{amssymb}
\usepackage[centertags]{amsmath}
\usepackage{amscd}
\usepackage{amsxtra}
\usepackage{latexsym}
\usepackage[shiftmargins]{vmargin}
\setpapersize{A4}
\setmarginsrb{35mm}{25mm}{30mm}{25mm}{10mm}{2\baselineskip}{10mm}{3\baselineskip}

\newtheorem{thm}{Theorem}[section]
\newtheorem{prop}[thm]{Proposition}
\newtheorem{lem}[thm]{Lemma}
\newtheorem{cor}[thm]{Corollary}
\newtheorem{definition}[thm]{Definition}
\newtheorem{remark}[thm]{Remark}
\newtheorem{example}[thm]{Example}

\newenvironment{rem}{\begin{remark}\begin{rm}}{\end{rm}\end{remark}}
\newenvironment{defi}{\begin{definition}\begin{rm}}{\end{rm}\end{definition}}

\newcommand{\bt}{\begin{thm}}
\newcommand{\br}{\begin{rem}}
\newcommand{\bc}{\begin{cor}}
\newcommand{\bl}{\begin{lem}}
\newcommand{\bd}{\begin{defi}}
\newcommand{\bp}{\begin{prop}}

\newcommand{\et}{\end{thm}}
\newcommand{\er}{\end{rem}}
\newcommand{\ec}{\end{cor}}
\newcommand{\el}{\end{lem}}
\newcommand{\ed}{\end{defi}}          
\newcommand{\ep}{\end{prop}}

\def\qedbox{\mbox{\large $\Box$}}
\newenvironment{Pf}{\noindent{\it
       Proof:$\;\; $}}{\rule{0mm}{0mm}\nolinebreak\hfill\qedbox \vspace{12pt}}
\newcommand{\bpf}{\begin{Pf}}
\newcommand{\epf}{\end{Pf}}

\newcommand{\boxma}{\\\begin{minipage}[b]{14cm}}
\newcommand{\boxme}{\end{minipage} \hfill $\Box$}
\newcommand{\mb}[1]{\mathbb{#1}}
\newcommand{\rn}{\ensuremath{\mb{R}^{n}}}
\newcommand{\hn}{\ensuremath{\mathcal{H}^{n}}}
\newcommand{\hne}{\ensuremath{\mathcal{H}^{n+1}}}

\newcommand{\R}{\ensuremath{\mb{R}}}

\newcommand{\rne}{{\ensuremath{\mb{R}^{n+1}}}}
\newcommand{\rnz}{{\ensuremath{\mb{R}^{n+2}}}}
\newcommand{\res}{\ensuremath{\,\textsf{\small \upshape L}\,}}
\newcommand{\p}{\partial}
\newcommand{\beq}{\begin{equation}{}}
\newcommand{\eeq}{\end{equation}{}}
\newcommand{\tr}{\, \triangle \,}

\newcommand{\eps}{\varepsilon}

\newcommand{\inl}{\int\limits}
\newcommand{\rhp}{\rightharpoonup}
\renewcommand{\leq}{\leqslant}
\renewcommand{\geq}{\geqslant}
\newcommand{\ra}{\rightarrow}
\newcommand{\nd}[1]{\frac{1}{#1}}

\newcommand{\dmu}{d\mu}
\newcommand{\cp}{\mathrm{Cap}}
\newcommand{\graph}{\mathrm{graph}\,}
\newcommand{\supp}{\mathrm{supp}\,}

\renewcommand{\div}{\mathrm{div}}
\newcommand{\lan}{\langle}
\newcommand{\ran}{\rangle}
\newcommand{\IR}{\mb{R}}
\newcommand{\CH}{\mathcal{H}}
\newcommand{\CL}{\mathcal{L}}
\newcommand{\tOmega}{{\tilde\Omega}}

\newcommand{\Net}{{N_t^{\eps_i}}}
\newcommand{\Het}{{\vec H_t^{\eps_i}}}
\newcommand{\dmuet}{\dmu_{\eps_i,t}}
\newcommand{\nuet}{{\nu_t^{\eps_i}}}
\newcommand{\Ue}{{U^{\eps_i}}}
\newcommand{\ale}{{\alpha^{\eps_i}}}
\newcommand{\tUet}{{\{t_1\leq \Ue \leq t_2\}}}
\newcommand{\tUt}{{\{t_1\leq U \leq t_2\}}}
\newcommand{\lUt}{{\{U=t\}}}
\newcommand{\tGmt}{{\tilde\Gamma_t}}
\newcommand{\tGmto}{{\tilde\Gamma_{t_1}}}
\newcommand{\tGmtt}{{\tilde\Gamma_{t_2}}}
 
\newcommand{\Gmt}{{\Gamma_t}}
\newcommand{\Gmto}{{\Gamma_{t_1}}}
\newcommand{\Gmtt}{{\Gamma_{t_2}}}

\newcommand{\Ht}{{\vec H_t}}
\newcommand{\tdmu}{d\tilde\mu}
\newcommand{\dmut}{\dmu_{t}}
\newcommand{\tdmut}{\tdmu_{t}}
\newcommand{\nut}{{\nu_t}}

\begin{document}

\title[No mass drop for MCF]{No mass drop for mean curvature flow of mean convex hypersurfaces}

\author{Jan Metzger}
\address{Jan Metzger, Max-Planck-Institut f\"ur Gravitationsphysik, 
Am M\"uhlenberg 1,\newline 14476 Golm, Germany}
\email{jan.metzger@aei.mpg.de}


\author{Felix Schulze}
\address{Felix Schulze, FU Berlin, Arnimallee 2-6, 
  14195 Berlin, Germany}
\email{felix.schulze@math.fu-berlin.de}




\begin{abstract}
A possible evolution of a compact hypersurface in $\rne$ by mean curvature
past singularities is defined via the level set flow. In the
case that the initial hypersurface has positive mean curvature, we show
that the Brakke flow associated to the level set flow is actually a
Brakke flow with equality. We obtain as a consequence that no mass drop
can occur along such a flow. As a further application of the
techniques used above we give a new variational formulation for mean
curvature flow of mean convex hypersurfaces. 
\end{abstract}

\maketitle
\section{Introduction}

Let $M\subset \rne$ be a smooth, compact $n$-dimensional submanifold without
boundary, and let $(M_t)_{t\in[0,T)}$ be the maximal smooth evolution of $M$
by mean curvature flow. Since $M$ is compact, the maximal time of existence
$T$ is finite, and in general the flow will develop singularities
before the surfaces vanish. One way to define a weak solution past singularities
is the level set flow of Chen-Giga-Goto \cite{ChenGigaGoto} and
Evans-Spruck \cite{EvansSpruckI}. Let us briefly recall one way of defining
the level set flow. It uses the so-called avoidance
principle: If two smooth mean curvature flows (where at least one of them is
compact) are disjoint at time $t_0$, then they remain so for all
times $t>t_0$. A weak mean curvature flow, generated by $M$, is a closed
subset $\mathcal{M}$ of space time $\rne\times \R^+$ such that for
$$M_t:=\{x\ |\ (x,t)\in \mathcal{M}\}$$
we have $M_0=M$, and the family of sets $(M_t)_{t\geq 0}$ satisfy the
above avoidance principle with respect to any smooth mean curvature flow. The
level set flow of $M$ is then characterized as the unique maximal weak
mean curvature flow generated by $M$. Assume now that $M$ has
nonnegative mean curvature. Following \cite{White00} and \cite{Ilmanen}, the level set
flow $\mathcal{M}$, generated by $M$, has further properties. Let 
$K\subset \rne \times \R^+$ be the compact set enclosed by the 
level set flow $\mathcal{M}$, such that $\partial K = \mathcal{M}$. 
Then $M_t= \partial K_t$, where $K_t=\{x \in \rne\ |\ (x,t)\in K\}$,
and the family of Radon measures
\beq \label{intro.1}
 \mu_t:=\hn\res \partial^*\!K_t 
\eeq
constitutes a Brakke flow. It has further regularity properties
for almost every $t$: $M_t=\partial^*\!K_t$ up to $\hn$-measure zero
and $M_t$ is a unit density, $n$-rectifiable varifold which carries a
weak mean curvature $\vec{H}$. The fact that $(\mu_t)_{t\geq 0}$ is a Brakke flow
can be characterized as follows. Given any $\phi \in
C_c^2(\rne;\R^+)$ the following inequality holds for every $t>0$
\beq\label{intro.2}
\bar{D}_t\mu_t(\phi)\leq \int -\phi |\vec{H}|^2 + \langle\nabla \phi,
\vec{H}\rangle \, d\mu_t
\eeq
where $\bar{D}_t$ denotes the upper derivative at time $t$ and we take
the left hand side to be $-\infty$, if $\mu_t$ is not $n$-rectifiable, or doesn't
carry a weak mean curvature. Note that in case $M_t$ moves smoothly by
mean curvature, $\bar{D}_t$ is just the usual derivative and we have
equality in (\ref{intro.2}). We can now state our main result.

\bt\label{intro.t1} Let $\Omega$ be an open and bounded subset of
$\rne$. Assume further that
 $M =\partial \Omega$ is a closed submanifold of $\rne$ of class
$C^{1}$, carrying a nonnegative weak mean curvature in $L^2$. 
Then the family of Radon measures $(\mu_t)_{t\geq 0}$
associated to the level set flow of $M$ is a Brakke flow with
equality in the sense that
\beq\label{intro.3}
\mu_{t_2}(\phi) - \mu_{t_1}(\phi) = \inl_{t_1}^{t_2}\!\!\!\int -\phi H^2 + \langle\nabla \phi,
H\rangle \, d\mu_t\, dt\ ,
\eeq 
for $0\leq t_1\leq t_2$ and any $\phi \in C_c^2(\rne)$. 
\et
This implies the following property, which is known as no mass drop.

\bc\label{intro.c1}
The family of Radon measures $(\mu_t)_{0\leq t\leq T}$, where $T$ is
the maximal time of existence of the level set flow, is continuous in
time. Furthermore $\lim_{t\nearrow T}\mu_t(\phi) = 0$ for any $\phi
\in C_c^2(\rne;\R^+)$. 
\ec

As a further application of the techniques used in the proof above, we give a variational
formulation for mean curvature flow of mean convex surfaces which is
similar to the variational principle applied by Huisken-Ilmanen to
define weak solutions to the inverse mean curvature flow, see
\cite{HuiskenIlmanen}.\\[1ex]
Since $\partial \Omega$ has nonnegative mean curvature and $\Omega$ is
compact, we expect that that by the strong maximum principle the 
mean curvature of the evolving surfaces becomes immediately strictly
positive and remains so. The surfaces can then be given as level-sets
of a continuous function $u:\bar{\Omega}\ra \R^+$, $u=0$ on $\partial
\Omega$ via 
$$\partial\{x\in\Omega\ |\ u(x)>t\}\ $$
and $u$ satisfies the degenerate elliptic equation
\beq\label{eq:weak_mcf}\tag{\ensuremath{\star}} 
\text{div}\Big(\frac{Du}{|Du|}\Big) = -\frac{1}{|Du|}\ .
\eeq
Note that if $u$ is smooth at a point $x \in \Omega$ with $Du(x)\neq
0$, this equation just states that the level sets of $u$ near $x$ are flowing smoothly by
mean curvature. To give this equation a variational structure we
proceed as follows. Given a function $w\in C^{0,1}(\bar{\Omega})$ such
that $|Dw|^{-1}\in L^1(\Omega)$ we define the functional
$$ J_w(v):= \int_\Omega |Dv| - \frac{v}{|Dw|}\, dx\ ,  $$
for any Lipschitz continuous function $v$ on $\Omega$ such that
$\{w\neq v\}\Subset \Omega$. We then say that such a function $w$ is a
weak solution to $(\star)$ on $\Omega$ if 
$$J_w(w)\leq J_w(v)$$
for any such $v$ as above and $w$ fulfills the boundary conditions
\beq\label{bc} w>0\ \ \text{on}\ \Omega, \ \ \{w=0\}=\p\Omega\ .\eeq
For an $\Omega$ satisfying the conditions of theorem \ref{intro.t1}
we show that the level set flow of $\p\Omega$ can be described as the
level-sets of a function $u$ satisfying $(\ref{bc})$ with
$|Du|^{-1}\in L^1(\Omega)$.

\bt\label{intro.t2} Let $\Omega$ be as in theorem \ref{intro.t1}. Then
the level set flow $u:\bar{\Omega}\ra \R$ is the unique weak solution
to $(\star)$ on $\Omega$.
\et

Mean curvature flow in the varifold setting was pioneered by Brakke
\cite{Brakke}. Aside from the classical PDE setting, see for example
\cite{GageHamilton86}, \cite{Huisken84}, the level set approach using
viscosity techniques proved to be fruitful, see
\cite{ChenGigaGoto}, \cite{EvansSpruckI}. The level set flow approach
has the advantage to define the evolution by mean curvature for
any closed subset of $\rne$. An alternative approach using Geometric
Measure Theory together with an approximative variational functional
to yield Brakke flow solutions is given in \cite{Ilmanen}. In \cite{EvansSpruckIV},
or equivalently \cite{Ilmanen}, a connection between the level set flow and the
varifold setting is drawn: For a given initial function $u_0$ on $\rne$ they show
that the level set flow of a.e. level set of $u_0$ constitutes a
Brakke flow. For a compact initial hypersurface that has nonnegative
mean curvature this gives that its level set flow constitutes a Brakke
flow. In our work here we partly use this techniques to show that
in this case one actually has equality in the Brakke flow definition.
 \\[2ex]
{\bf Outline}. In \S\ref{pr} we recall a way 
of defining the level set flow by Evans-Spruck for initial
hypersurfaces with positive mean curvature. We work out some geometric
consequences of approximation by elliptic regularisation. The
approximating solutions $u^\varepsilon$ have the important geometric
property that, scaled appropriately, they constitute a smooth,
graphical, translating solution to mean curvature flow in
$\Omega\times\R$. Writing these translating graphs again as level sets
of a function $U^\varepsilon$ on $\Omega\times\R$ this yields an
approximation of the level set flow $U$, where $U$ is the constant
extension of $u$ in the $e_{n+2}$-direction, by smooth level set flows.\\[1ex]
In \S \ref{arg} we show that the obstacle to $U$ being a Brakke
flow with equality can be characterized by the possible existence of a
nonnegative radon measure $\gamma$ such that for a subsequence $\varepsilon_i \ra 0$,
\beq\label{intro.4}
|DU^{\epsilon_i}|^{-1}\rightharpoonup |DU|^{-1}+\gamma 
\eeq
on $\Omega\times \R$. 
Using the regularity results of White,
\cite{White00}, we can furthermore deduce that the $1$-capacity of the
support of $\gamma$ has to vanish. Since the limit flow is still a kind of
Brakke flow with equality, where this incorporates the defect measure
$\gamma$, we can apply this limit equation to show that $\gamma$
actually has to vanish entirely.\\[1ex]
To be able to state theorem
\ref{intro.t1} not only for boundaries $\partial\Omega$ which are smooth
with positive mean curvature, we show that the
level set flow of any boundary $\partial\Omega$ as in theorem \ref{intro.t1} is actually
smooth for positive times close enough to zero and has positive mean curvature.
\\[1ex]
In \S \ref{varprinc} we employ that $(\ref{intro.4})$ holds with
$\gamma \equiv 0$ and use the approximation by smooth level set flows
one dimension higher to
show that $u$ also is a weak solution to $(\star)$. Since the
variational principle implies that a weak solution constitutes a Brakke
flow with equality, we can apply the avoidance principle for Brakke
flows to show uniqueness of such a weak solution.


\section{Preliminaries}\label{pr}

In the case that our initial hypersurface $M$ is the smooth boundary of an 
open and bounded set $\Omega\subset \rne$, and $M$ has mean curvature
$H>0$, the level set flow $\mathcal{M}$ generated by $M$ can be written
as the graph of a continuous function $u:\Omega \ra \R^+$. In
\cite{EvansSpruckI} it is shown that $u$ is the unique continuous
viscosity solution of
\beq\label{pr.1}\begin{split}
\qquad\bigg(\delta^{ij} -\frac{D^iu D^ju}{|Du|^2}\bigg)D_iD_ju =&\ -1\qquad
\text{in}\ \Omega\ ,\\
\qquad u =&\ 0 \qquad \ \ \ \ \text{on}\ \partial\Omega\ .
\end{split}\eeq
Note that if $u$ is smooth at a point $x \in \Omega$ with $Du(x)\neq
0$, equation (\ref{pr.1}) is identical to (\ref{eq:weak_mcf}).
To prove the existence of a solution to (\ref{pr.1})
Evans and Spruck use the method of elliptic regularisation. Since
$(\star)$ is degenerate everywhere on $\Omega$ and singular at points 
where $Du = 0$, one replaces it 
by the following nondegenerate PDE,
\beq\label{eq:eps_mcf}\tag{\ensuremath{\star_\eps}} \begin{split}
\text{div}\bigg(\frac{Du^\eps}{\sqrt{\eps^2+|Du^\eps|^2}}\bigg) =&\
-\frac{1}{\sqrt{\eps^2+|Du^\eps|^2}} \qquad \text{in}\ \Omega\ ,\\
u^\eps =&\ 0 \qquad \qquad \qquad \qquad \ \ \ \text{on}\ \partial \Omega\ ,
\end{split}\eeq
for some small $\eps>0$. Since the mean curvature of the boundary
$\partial \Omega$ is strictly positive one can construct barriers at
the boundary, which yield, together with a maximum principle for the
gradient, uniform a-priori gradient bounds for solutions $u^\eps$, provided
$0<\eps<1$. Applying De-Giorgi/Nash/Moser and Schauder estimates
together with a continuity argument then gives existence of solutions
to $(\star_\eps)$ for small $\eps>0$. Furthermore, by Arzela-Ascoli, there is a
sequence $\eps_i\ra 0$ such that $u^{\eps_i}\ra \tilde{u}$ and
$\tilde{u}\in C^{0,1}(\Omega)$ is a solution to $(\ref{pr.1})$. By
uniqueness $\tilde{u}=u$ and $\lim_{\eps\ra 0} u^{\eps} = u$. Aside
from the $\eps$-independent gradient estimate there is also a uniform integral
estimate for the RHS of $(\star_\eps)$:

\bl\label{pr.l1} For any solution $u^\eps$ of $(\star_\eps)$ we have the
bound
\beq\label{pr.2}
\inl_\Omega \frac{1}{\sqrt{\eps^2+|Du^\eps|^2}} \, dx \leq |\partial
\Omega|
\eeq 
\el
\bpf Choose a smooth function $\varphi$, $0\leq \varphi \leq 1$, such that $\varphi=1$
on $\Omega_\delta:=\{x\in \Omega\ | \ \text{dist}(x,\partial \Omega)>
\delta\}$, $\varphi=0$ on $\partial \Omega$ and $|D\varphi|\leq \gamma/\delta$ for
some $\gamma >1, \delta > 0$. Multiplying $(\star_\eps)$ with $\varphi$ and
integrating by parts we obtain that
$$ \inl_{\Omega_\delta} \frac{1}{\sqrt{\eps^2+|Du^\eps|^2}} \, dx \leq
\inl_{\Omega\setminus\Omega_\delta} |D\varphi|
\frac{|Du^\eps|}{\sqrt{\eps^2+|Du^\eps|^2}}\, dx \leq
\frac{\gamma}{\delta}|\Omega\setminus\Omega_\delta| \ .$$
Then taking the limit $\gamma \ra 1$ and $\delta \ra 0$ gives the
estimate.
\epf

The sequence $u^{\eps_i}$ is bounded in
$C^1(\overline{\Omega})$ and converges uniformly to $u$, which thus is in
$C^{0,1}(\overline{\Omega})$. Thus also
$$Du^{\eps_i} \rhp Du \qquad \text{weakly}-*\ \ \text{in}\
L^{\infty}(\Omega;\rne)\ .$$ 
Together with the estimate $(\ref{pr.2})$ we can then apply Theorem 3.1 in
\cite{EvansSpruckIV}, \S 3 (by a slight
modification of the proof there, since here the functions $u^{\eps}$
are not defined on all of $\rne$) to obtain:
 
\bp[Evans/Spruck]\label{pr.p1}
We have the following convergence:
\beq\label{pr.3}
|Du^{\eps_i}| \rhp |Du|\qquad \text{weakly}-*\ \text{in}\
L^{\infty}(\Omega)\ ,
\eeq
\beq\label{pr.4}
\frac{Du^{\eps_i}}{\sqrt{\eps_i^2+|Du^{\eps_i}|^2}} \ra
\frac{Du}{|Du|}\ \ \text{in}\ L^p(\{|Du|\neq 0\}\cap \Omega ; \rne\})
\eeq
for any $ p\geq 1$ .
\ep
The following Lemma and Proposition are a direct consequence of this
strengthened convergence. In fact they are a variation of Lemma 4.2
and step 4 in the proof of Theorem 5.2 in \cite{EvansSpruckIV}. 
\bl\label{pr.l2}
$$\hne(\{x\in\Omega\ |\ Du(x)=0\})=0\ . $$
\el
\bpf
Let $A:=\{Du=0\}\subset \Omega$. By Lemma \ref{pr.l1} we have
\begin{equation*}\begin{split}
\inl_{A} 1\, dx \ =& \ \lim_{\eps_i\rightarrow
  0}\inl_{A}\sqrt{\eps_i^2+|Du^{\eps_i}|^2}^{\ -\frac{1}{2}}
\sqrt{\eps_i^2+|Du^{\eps_i}|^2}^\frac{1}{2}\ dx\\
\leq & \ \limsup_{\eps_i\rightarrow
  0}\bigg(\inl_{A}\sqrt{\eps_i^2+|Du^{\eps_i}|^2}^{\ -1}\,
  dx\bigg)^\frac{1}{2}
\cdot \bigg(\inl_{A}\sqrt{\eps_i^2+|Du^{\eps_i}|^2}\, dx\bigg)^\frac{1}{2}\\
\leq& \ C \ \limsup_{\eps_i\rightarrow
  0}\bigg(\inl_{A}\sqrt{\eps_i^2+|Du^{\eps_i}|^2}\ dx\bigg)^\frac{1}{2}\ .
\end{split}\end{equation*}
Since $\big|\sqrt{\eps_i^2+|Du^{\eps_i}|^2}-|Du^{\eps_i}|\big|\leq
\eps_i$ we have by (\ref{pr.3}) that
$$\hne(A)\leq C \limsup_{\eps_i\rightarrow 0}\bigg(
\inl_{A}|Du^{\eps_i}|\ dx \bigg)^\frac{1}{2}=0\ .$$
\epf

Thus, by (\ref{pr.4}), 
\beq\label{pr.4a} \frac{Du^{\eps_i}}{\sqrt{\eps_i^2+|Du^{\eps_i}|^2}} \ra
\frac{Du}{|Du|}
\eeq in $ L^p(\Omega)$ for any $p\geq 1$.

\bp\label{pr.p2} For $\phi \in L^\infty(\Omega),\ \phi\geq 0$ we have:
$$ \inl_{\Omega}\frac{\phi}{|Du|}\, dx\ \leq \ \liminf_{\eps_i\rightarrow 0}
\inl_{\Omega}\frac{\phi}{\sqrt{\eps_i^2+|Du^{\eps_i}|^2}}\, dx\ .$$
\ep
\bpf Let $\phi,\psi \in L^\infty(\Omega),\ \phi,\psi\geq 0$. One obtains 
\begin{align*}
\inl_{\Omega}\phi\psi&\, dx =\ \lim_{\eps_i\rightarrow
0}\inl_{\Omega}\Big(\phi^\nd{2}\psi\sqrt{\eps_i^2+|Du^{\eps_i}|^2}^\nd{2}\Big)\Big(
\phi^\frac{1}{2}\sqrt{\eps_i^2+|Du^{\eps_i}|^2}^{\ -\nd{2}}\,
dx\Big)\\ \displaybreak[0]
\leq&\ \liminf_{\eps_i\rightarrow
0}\bigg(\inl_{\Omega}\phi\psi^{2}\sqrt{\eps_i^2+|Du^{\eps_i}|^2}\ dx\bigg)^\nd{2}
\bigg(\inl_{\Omega}\phi\sqrt{\eps_i^2+|Du^{\eps_i}|^2}^{\ -1}\
dx\bigg)^\frac{1}{2}\\ \displaybreak[0]
=&\ \bigg(\inl_{\Omega}\phi\psi^{2}|Du|\ dx\bigg)^\nd{2} \liminf_{\eps_i\rightarrow
0}\bigg(\inl_{\Omega}\phi\sqrt{\eps_i^2+|Du^{\eps_i}|^2}^{\ -1}\ dx\bigg)^\frac{1}{2}\ .
\end{align*}
Now choose $\psi:=\varphi_m\big(|Du|^{-1}\big)$ with $\varphi_m :\R\rightarrow \R$,
$$ \varphi_m(z)=\begin{cases}\  m\ \ &\ \ \ \ \  z\geq m\\\ z\ \ &-m\leq z\leq m\\
\ -m\ \ &\ \ \ \ \  z\leq -m\ \ . \end{cases}$$
Since $\psi \leq |Du|^{-1}$ we obtain by the calculation above
$$\bigg(\inl_{\Omega}\phi\psi\ dx\bigg)^\frac{1}{2}\leq \liminf_{\eps_i\rightarrow
0}\bigg(\inl_{\Omega}\phi\sqrt{\eps_i^2+|Du^{\eps_i}|^2}^{\ -1}\, dx\bigg)^\frac{1}{2}\ ,$$
which by the monotone convergence theorem for $m\rightarrow \infty$
proves the claim. \epf

The equation $(\star_\eps)$ also has a geometric interpretation. It
implies that the downward translating graphs
$$N_t^\eps:=\text{graph}\bigg(\frac{u^\eps(x)}{\eps}-\frac{t}{\eps}\bigg),\qquad
-\infty<t<\infty $$
are hypersurfaces in $\Omega\times \R$ flowing smoothly by mean
curvature flow. To verify
this, define the function 
\beq\label{pr.5}U^\eps(x,z):=u^\eps(x)-\eps z,\qquad (x,z)\in
\Omega\times \R\ ,
\eeq
such that $\{U^\eps=t\}=N^\eps_t$. It is easily checked that
$U^\eps$ satisfies $(\star)$ on $\Omega\times\R$ if and only if
$u^\eps$ satisfies $(\star_\eps)$ on $\Omega$. Since the level set
flow $u$ is non-fattening (see \cite{Ilmanen}), we have that
$$ K_t=\{x\in\Omega\ |\ u(x)\geq t \} $$
and 
$$ \hne(\{u=t\}) = 0 $$
for all $t\geq 0$, which implies that
$$\p^*K_t = \p^*\{u>t\} \ \qquad \forall\ t\geq 0\ .$$ Note that $u$
is Lipschitz continuous and thus also a $BV$-function. By comparing
the coarea formula for Lipschitz functions and for BV-functions we see
that 
$$ \p^*\{u>t\} = \{u=t\}$$
up to $\hn$-measure zero for almost every $t$. Then define the
family $\tilde{\mu}_t$ of $(n+1)$-rectifiable radon measures on
$\Omega \times \R$ by
$$ \tilde{\mu}_t:=\hne\res (\p^*K_t\times \R) \ ,$$ 
i.e. $\tilde{\mu}_t=\mu_t\otimes \mathcal{L}^1$, where $\mu_t$ is the
Brakke flow associated to the level set flow $u$. We can now make
precise in what sense the translating graphs $N^\eps_t$ approximate
the Brakke flow $\mu_t$. For a proof of the following Proposition see
\S 5 in \cite{Iso06}.

\bp\label{pr.p3} Let $\eps_i\ra 0$. Then for almost all $t \geq 0$ we have  
\beq\label{pr.6}
\hne\res N^{\eps_i}_t \ra \tilde{\mu}_t
\eeq 
in the sense of radon measures. Even more for almost every $t\geq 0$
there is a subsequence $\{\eps_j\}$, depending on $t$, such that
\beq\label{pr.7}
N^{\eps_j}_t\ra \p^*K_t\times \R
\eeq
in the sense of varifolds. Furthermore the rectifiable sets $\p^*K_t$,
seen as unit density $n$-rectifiable varifolds carry for a.e. $t\in
[0,T)$ a weak mean curvature $\vec{H}_t \in L^2(\p^*K_t, \hn)$.
\ep


\section{The argument}\label{arg}
Throughout this section we will work with the level set flows of $\Ue$
in $\tOmega:=\Omega\times[0,1]$. We denote $U(x,z) = u(x)$.
Furthermore, let $\nut$ and $\nuet$ be the normal vectors to the
level set flows $U$ and $\Ue$ respectively and  denote
$\Gamma_t:=\p^*K_t$ and $\tGmt:=\Gamma_t\times\IR$. From lemma \ref{pr.l1},
proposition \ref{pr.p1}, and lemma \ref{pr.l2}, we derive the
following facts.
\begin{lem}
  \begin{gather}
    \label{eq:bdd}
    \int_{\tOmega} |D\Ue|^{-1} dx \leq |\p\Omega|,
    \\
    \label{eq:weak_l1}
    |D\Ue| \rhp |DU|\qquad\text{weakly-}*\text{\ in\ }
    L^\infty(\tOmega), \text{and}
    \\
    \label{eq:nuconv}
    \nuet\to\nut\qquad\text{in\ }L^1(\tOmega,\IR^{n+2}).
  \end{gather}
\end{lem}
We now introduce several Radon measures on $\tOmega$ which will be
central in the proof of theorem~\ref{intro.t1}. First, denote
\begin{equation*}
  \begin{split}
    \ale : = |D\Ue|^{-1}\,dx^{n+2},\quad \text{and}
    \\
    \alpha := |DU|^{-1}\,dx^{n+2}.
  \end{split}
\end{equation*}
Since for all $K\subset\subset\tOmega$
\begin{equation*}
  \ale(K) = \int_K |D\Ue|^{-1} dx  \leq  C(K) 
\end{equation*}
by equation~\eqref{eq:bdd}, we know that the $\ale$ have a convergent
subsequence, i.~e. we can assume 
\begin{equation*}
  \ale \rhp \beta 
\end{equation*}
in the sense of Radon measures. We will clarify the relation of
$\alpha$ and $\beta$ subsequently. Note that in view of
proposition~\ref{pr.p2} we find that $\beta\geq\alpha$. We will denote
the defect measure, the difference of $\alpha$ and $\beta$, by
\begin{equation*}
  \gamma = \beta-\alpha,
\end{equation*}
which is a non-negative Radon measure. 

Before attempting the proof of theorem \ref{intro.t1}, we collect some
observations about the level set flow $U$. 
\begin{lem}
  \label{lem:H}
  Let $\Ht$ denote the mean curvature vector of the level set flow
  $U$. Then for all smooth vector fields $X$ with
  compact support in $\tOmega$,
  \begin{equation*}
    \int_\tOmega\lan\Ht,X\ran |DU|\,dx
    = \int_\tOmega
    \lan\tfrac{DU}{|DU|},X\ran\,dx,
  \end{equation*}
  i.e. $\Ht$ agrees with $\frac{DU}{|DU|^2}$ almost everywhere, with
  respect to the measure $|DU|\,dx$.
\end{lem}
\begin{proof}
  Let $X$ be a smooth vector field with compact support in
  $\tOmega$. Then for all $\eps_i>0$ there exists $T>0$,
  such that
  \begin{equation*}
    \supp(X)\cap \Net = \emptyset\quad\forall t\not\in[-T,T].
  \end{equation*}
  Denote by $\nuet$ and $\Het$ the downward normal and mean
  curvature vector of $\Net$. Using the fact that the $\Net$ are
  smooth and the co-area formula for the level-set function $\Ue$, we
  compute
  \begin{equation}
    \label{eq:Hlemma1}
    \begin{split}
      \int_{-\infty}^{\infty} \int_\Net \lan\Het, X\ran
      \dmuet\, dt
      &=
      -\int_{-T}^{T} \int_\Net \div_\Net(X)\dmuet\, dt
      \\
      &=
      -\int_{\tOmega} \big(\div_\rnz X - \lan D_\nuet X,
      \nuet\ran\big)|D\Ue|\, dx.
    \end{split}
  \end{equation}
  We claim that
  \begin{equation}
    \label{eq:Hlemma2}
    \int_{\tOmega} \lan D_\nuet X,\nuet\ran|D\Ue|\, dx
    \to
    \int_{\tOmega}\lan D_\nut X,\nut\ran |DU|\, dx,
  \end{equation}
  as $i\to\infty$. Indeed, we can write 
  \begin{equation*}
    \begin{split}
      &
      \int_{\tOmega} \lan D_\nuet X,\nuet\ran |D\Ue|\,dx
      \\
      &\quad=
      \int_{\tOmega} \lan D_{\nuet-\nut} X,\nuet\ran |D\Ue|
      + \lan D_\nut X,\nuet-\nut\ran |D\Ue|\,dx
      + \int_{\tOmega}\lan D_\nut X,\nut\ran |D\Ue|\,dx.
    \end{split}
  \end{equation*}
  From the fact that $\lan D_\nut X,\nut\ran\in
  L^1(\tOmega)$ and $|D\Ue|\rhp|DU|$ weakly-$*$ in
  $L^\infty(\tOmega)$, we infer that the last term in the above equation
  converges to $\int_{\tOmega} \lan D_\nut X,\nut\ran
  |DU|dx$. The first two terms go to zero, since $\nuet\to \nut$ in
  $L^1(\tOmega)$ and the respective factors are bounded. Thus we
  established \eqref{eq:Hlemma2}. Subsequently, we will use the
  co-area-formula in the form
  \begin{equation*}
    \int_{\IR^{n+2}} f\,dx = \int_\IR \int_{\{U=t\}} f |DU|^{-1}d\CH^{n+1}\,dt
  \end{equation*}
  for $f\in L^{\infty}$. This is justified by \cite[Theorem 2,
  p.117]{Evans-Gariepy}, since in view of lemma~\ref{pr.l1} and
  proposition~\ref{pr.p2} we have that $\int_{\tOmega}
  |DU|^{-1} dx <\infty$. Thus we can compute
  \begin{equation}
    \label{eq:Hlemma3}
    \begin{split}      
      -\int_{\tOmega}\big(\div_\rnz X - \lan D_\nut X,
      \nut\ran\big) |DU|dx
      &=
      -\int_{\tOmega} \div_{\tGmt}X |DU| dx
      \\
      &= - \int_{\{t>0\}} \int_{\tGmt}\div_{\tGmt} X \tdmu_t\,dt
      \\
      &=
      \int_{\{t>0\}} \int_{\tGmt} \lan \Ht,X\ran \tdmu_t\,dt.
    \end{split}
  \end{equation}
  On the other hand, since the $\Net$ constitute a smooth level set
  flow, we have that
  \begin{equation*}
    \begin{split}
      \int_{-T}^T\int_\Net \lan \Het, X\ran \dmuet dt
      &=
      \int_{-T}^T\int_\Net \lan \nuet,X\ran |D\Ue|^{-1}\dmuet\,dt
      \\
      &=
      \int_{\tOmega} \lan\nuet,X\ran dx      
    \end{split}
  \end{equation*}
  As $\nuet\to\nut$ in $L^1$, we thus find
  \begin{equation*}    
    \begin{split}
      \int_{-T}^T\int_\Net \lan \Het, X\ran \dmuet\, dt
      &\to \int_{\tOmega} \lan \nu_t, X\ran dx
      \\
      &= \int_{\{t>0\}} \int_{\tGmt}
      \lan\nu,X\ran|DU|^{-1}\tdmu_t\,dt.
    \end{split}
  \end{equation*}
  Combining this equation with \eqref{eq:Hlemma1},~\eqref{eq:Hlemma2},
  and~\eqref{eq:Hlemma3}, we find that for all $X$
  \begin{equation*}
    \int_{\{t>0\}} \int_{\tGmt} \lan \Ht,X\ran \tdmu_t\,dt
    =
    \int_{\{t>0\}} \int_{\tGmt} \lan\nut,X\ran|DU|^ {-1}\tdmu_t\,dt.
  \end{equation*}
  An application of the co-area formula yields the claimed identity.
\end{proof}
We are now set up to perform the central computation. Fix $\phi\in
C_c^\infty(\tOmega)$ and let $0<t_1<t_2$. We will adopt the
convention, that if $t>\sup_\Omega u$, then $\tGmt=\emptyset$. Note that
also $\Net\cap\tOmega=\emptyset$ if $t\not\in[-T,T]$, provided $T$ is
large enough.

Since $\Net$ is a smooth level set flow, we know that
\begin{equation}
  \label{eq:area_change1}
  \begin{split}
    &\int_{N^{\eps_i}_{t_2}} \phi\,\dmu_{\eps_i,t_2}
    -
    \int_{N^{\eps_i}_{t_1}} \phi\,\dmu_{\eps_i,t_1}
    \\
    &\quad=
    \int_{t_1}^{t_2} \int_\Net
    \lan \nabla\phi,\Het\ran - \phi |\Het|^2
    \,\dmuet\,dt
    \\
    &\quad=
    - \int_{t_1}^{t_2} \int_\Net
    \div_\Net(\nabla\phi) - \phi |\Het|^2
    \,\dmuet\,dt
    \\
    &\quad=
    - \int_{\tOmega\cap\tUet}
    \big(\div_\rnz(\nabla\phi) - \lan D_\nuet\nabla\phi,\nuet\ran\big)
    |D\Ue|\,dx
    \\
    &\qquad 
    - \int_{\tOmega\cap\tUet}
    \phi |D\Ue|^{-1}\,dx.
  \end{split}
\end{equation}
To take this computation to the limit as $i\to\infty$, observe that
non-fattening implies
\begin{equation*}
  \chi_{\tUet}\to \chi_\tUt\qquad\text{in}\quad L^1(\tOmega),
\end{equation*}
and thus the first term on the right hand side converges. Furthermore, for every $\delta>0$ consider the open set 
\begin{equation*}
  S_\delta = \{U \in (t_1-\delta,t_1+\delta)\}\cup\{U \in (t_2-\delta,t_2+\delta)\}.
\end{equation*}
As $\beta(\lUt) = 0$ for a.e. $t$, also for a.e. $t_1$,$t_2$:
\begin{equation*}
  \beta(\{U=t_1\}\cup\{U=t_2\}) = 0.
\end{equation*}
Since $\beta$ is a Radon measure, $\lim_{\delta\to 0} \beta(S_\delta)
= 0$. Hence, for every $\eta>0$ there exists $\delta>0$ such that
\begin{equation*}
  \beta(S_\delta) < \eta/2.
\end{equation*}
Therefore, there exists $N$ such that
\begin{equation*}
  \int_{S_\delta} |D\Ue|^{-1}\,dx\leq \eta \qquad\text{for all}\quad i\geq N.
\end{equation*}
In other words, as $i\to\infty$, eventually
$\ale(S_\delta)\leq\eta$. Thus
\begin{equation*}  
  \begin{split}
    &\left| \int_\tOmega\phi\chi_\tUet\,d\ale
      - \int_\tOmega \phi\chi_\tUt\,d\beta\right|
    \\
    &\quad
    \leq\left|\int_\tOmega\phi\chi_\tUet(1-\chi_{S_\delta})\,d\ale
      - \int_\tOmega \phi\chi_\tUt(1-\chi_{S_\delta})\,d\beta\right|
    \\
    &\qquad
    + 
    \left| \int_\tOmega\phi\chi_\tUet\chi_{S_\delta}\,d\ale
      - \int_\tOmega \phi\chi_\tUt\chi_{S_\delta}\,d\beta\right|.
  \end{split}
\end{equation*}
As $\Ue\to U$ uniformly, if $i$ is big enough,
\begin{equation*}
  \chi_\tUet(1-\chi_{S_\delta}) = \chi_\tUt(1-\chi_{S_\delta}),
\end{equation*}
which implies that the first term in the above computation goes to
zero, in view of the definition of $\beta$. The second term can be
estimated as follows
\begin{equation*}
  \begin{split}
    &\left| \int_\tOmega\phi\chi_\tUet\chi_{S_\delta}\,d\ale
      - \int_\tOmega \phi\chi_\tUt\chi_{S_\delta}\,d\beta\right|
    \\
    &\quad \leq \max|\phi|\big(\ale(S_\delta)+\beta(S_\delta)\big)
    \leq 2\eta\max|\phi|.
  \end{split}
\end{equation*}
In combination we find that 
\begin{equation*}
  \left| \int_\tOmega\phi\chi_\tUet\,d\ale
      - \int_\tOmega \phi\chi_\tUt\,d\beta\right| \to 0
\end{equation*}
as $i\to\infty$. In virtue of proposition~\ref{pr.p3} we can assume
that $N_{t_j}^{\eps_i} \to \tilde\Gamma_{t_j}$ for $j=1,2$ in the sense
of Radon measures. Then the above reasoning shows that
equation~\eqref{eq:area_change1} implies that
\begin{equation*}
  \begin{split}
    &\int_{\tGmtt} \phi\,\tdmu_{t_2}
    -
    \int_{\tGmto} \phi\,\tdmu_{t_1}
    \\
    &\quad=
    - \int_{\tOmega\cap\tUt} \big(\div_\rnz(\nabla\phi) - \lan
    D_\nut\nabla\phi,\nut\ran \big) |DU|\, dx - \int_{\tOmega\cap\tUt} \phi\,d\beta
    \\
    &\quad=
    \int_{\tOmega\cap\tUt}
    \lan\nabla\phi,\Ht\ran|DU|\,dx - \int_{\tOmega\cap\tUt}
    \phi\,d\beta.
  \end{split}
\end{equation*}
In view of lemma~\ref{lem:H} and the definition of the defect measure
$\gamma$, this yields
\begin{equation}
  \label{eq:area_change2}
  \begin{split}
    &\int_{\tGmtt} \phi\,\tdmu_{t_2}
    -
    \int_{\tGmto} \phi\,\tdmu_{t_1}
    \\
    &\quad=
    \int_{\tOmega\cap\tUt}
    \lan\nabla\phi,\tfrac{DU}{|DU|}\ran\,dx
    - \int_{\tOmega\cap\tUt}\phi\,d\alpha
    - \int_{\tOmega\cap\tUt}\phi\,d\gamma.
  \end{split}
\end{equation}
We now want to argue that the support of the defect measure $\gamma$
is very small. To do this, we introduce the notion of capacity
(cf.~\cite{Evans-Gariepy}).  For a set $A\subset\IR^n$, the
1-capacity, $\cp_1(A)$, is defined as follows:
\begin{equation*}
  \cp_1(A) = \inf \left\{ \int_{\IR^n} |Df|\,dx : f\geq 0, f\in C_c^\infty,
  A\subset \{f\geq 1\}^\circ \right\}.
\end{equation*}
It turns out that the 1-capacity of $\supp(\gamma)$ vanishes:
\begin{lem}
  \label{lem:cap}
  \begin{equation*}
    \cp_1(\supp(\gamma)\cap \tilde{\Omega}) = 0.
  \end{equation*}
\end{lem}
\begin{proof}
  From \cite{White00} we know that there is a closed singular set
  $\tilde S\subset\graph u$ of parabolic Hausdorff dimension at most
  $n-1$, such that outside of $\tilde S$ the sets $\{u=t\}$
  constitute a smooth level set flow. From the dimensionality we know
  that $\CH_\text{par}^n(\tilde S)=0$, where $\CH^n_\text{par}$
  denotes the $n$-dimensional parabolic Hausdorff measure. If we let
  $S := \Pi(\tilde S)$ be the projection of $\tilde S\subset\graph u$
  to $\Omega$, then we find that $\CH^n(S)=0$. In particular, $S$ is
  closed.

  Let $x_0\in \Omega\setminus S$. Then there exists a
  neighborhood $B = B_\delta(x_0)$ of $x_0$ such that $\graph u\big|_{B}$ is
  a smooth mean curvature flow. Thus, by Brakke's regularity theorem
  (cf. \cite{Ilmanen}), the $\Net\big|_{B\times\IR}$ converge smoothly
  on compact sub sets to $\tGmt\cap B\times\IR$.

  For $\phi\in C_c^\infty(B\times[0,1])$, we therefore conclude
  that
  \begin{equation*}    
    \begin{split}
      \int_{\tOmega} \phi|D\Ue|^{-1} dx
      &=
      \int_\IR\int_{\Net\cap B\times[0,1]} H^2 \phi\,\dmuet\,dt
      \\
      &\to \int_\IR\int_{\tGmt\cap B\times[0,1]}
      H^2\phi\,\tdmut\,dt
      =\int_{\Omega\times[0,1]}\phi|DU|^{-1}dx,
    \end{split}
  \end{equation*}
  as $i \to \infty$. Thus $x_0\not\in\supp(\gamma)$, which yields
  $\supp(\gamma)\subset S \times [0,1]$. Since $\CH^n(S)=0$ we find
  that $\CH^{n+1}(S\times\IR)=0$. Hence \cite[Section 5.6.3, Theorem
  3]{Evans-Gariepy} implies that $\cp_1(S\times[0,1])=0$.
\end{proof}
\begin{lem}
  \label{lem:ac1}
  For a.e. $0< t_1 < t_2$ and any $\phi\in C_c^\infty(\Omega\times\IR)$ we
  have
  \begin{equation*}
    \int_{\tGmtt}\phi\,\tdmu_{t_2}
    -
    \int_{\tGmto} \phi\,\tdmu_{t_1}
    =
    \int_{t_1}^{t_2}\int_{\tGmt} \lan\nabla\phi,\Ht\ran\tdmut\,dt
    - \int_{t_1}^{t_2}\int_{\tGmt}\phi |\Ht|^2 \tdmut\,dt.
  \end{equation*}
\end{lem}
\begin{proof}
  Without loss of generality we can assume $\supp(\phi)
  \subset\Omega\times[0,1]$.  Let $S = \supp(\gamma)$. Since
  $\cp_1(S)=0$ by lemma~\ref{lem:cap}, we can find functions $\eta_k\in
  C^\infty(\Omega\times\IR)$, such that $S\subset \{\eta_k\geq
  1\}^\circ$ and $\|\eta_k\|_{W^{1,1}(\rnz)} \to 0$ as
  $k\to\infty$.
  We can assume that the functions $\eta_k$ converge $\CL^{n+2}$ a.e.
  to $\eta\equiv 0$.
  
  Replace $\phi$ by  $(1-\eta_k)\phi$ in
  equation~\eqref{eq:area_change2}. Since $(1-\eta_k)=0$ on $S$, the
  term containing the defect measure $\gamma$ drops out, and we
  conclude
  \begin{equation}
    \label{eq:area_change3}
    \begin{split}
      &\int_{\tGmtt} (1-\eta_k)\phi\,\tdmu_{t_2}
      -
      \int_{\tGmto} (1-\eta_k)\phi\,\tdmu_{t_1}
      \\
      &\quad=
      \int_{\tOmega \cap \tUt}
      \big\lan \nabla  \big((1-\eta_k)\phi\big), \tfrac{DU}{|DU|} \big\ran\,dx
      - \int_{\tOmega\cap\tUt} (1-\eta_k)\phi\, d\alpha
      \\
      &\quad=
      \int_{\tOmega \cap \tUt}
      (1-\eta_k) \big\lan \nabla  \phi,  \tfrac{DU}{|DU|} \big\ran
      \,dx
      \\
      &\qquad
      - \int_{\tOmega \cap \tUt} (1-\eta_k) |DU|^{-1}\phi\,dx
      - \int_{\tOmega \cap \tUt} \phi \big\lan\nabla\eta_k, \tfrac{DU}{|DU|} \ran\,dx.
    \end{split}
  \end{equation}
  As $|DU|$ is bounded and $\eta_k\to 0$ in $L^1(\tOmega)$ we find
  \begin{equation*}
    \int_{\{t>0\}} \int_{\tGmt} |\eta_k \phi|\,\tdmut\,dt= \int_\tOmega
    |\eta_k\phi|\,|DU|\,dx \to 0,
  \end{equation*}
  as $k\to\infty$. Hence
  \begin{equation*}
    \int_{\tGmt} \eta_k\phi\,\tdmut \to 0\qquad\text{for a.e.}\quad t.
  \end{equation*}
  Thus for a.e. $0<t_1<t_2$, as $t\to\infty$ the left hand side of
  \eqref{eq:area_change3} converges to 
  \begin{equation*}
    \int_{\tGmtt} (1-\eta_k)\phi\,\tdmu_{t_2}
    -
    \int_{\tGmto} (1-\eta_k)\phi\,\tdmu_{t_1}
    \to
    \int_{\tGmtt} \phi\,\tdmu_{t_2}
    -
    \int_{\tGmto} \phi\,\tdmu_{t_1}.
  \end{equation*}
  To deal with the right hand side of~\eqref{eq:area_change3}, note
  that as $\eta_k\to 0$ a.e. and $\phi|DU|^{-1}$ is integrable, the
  second integrand converges  
  \begin{equation*}
    \int_{\tOmega \cap \tUt} (1-\eta_k) |DU|^{-1}\phi\,dx \to
    \int_{\tOmega \cap \tUt} |DU|^{-1}\phi\,dx
  \end{equation*}
  in view of the dominated convergence theorem. The other integrands
  converge in view of $\eta_k\to 0$ in $W^{1,1}(\rnz)$.

  Therefore in the limit as $k\to\infty$, equation \eqref{eq:area_change3}
  turns into the claimed identity.
\end{proof}
As a corollary of the proof of the previous lemma we find that the
defect measure $\gamma$ is in fact zero and we have convergence
$\ale\to\alpha$:
\begin{cor}\label{cor:measureconv}
  $|D\Ue|^{-1}dx \to |DU|^{-1}dx$ in the sense of Radon measures.
\end{cor}
The next lemma removes the extra dimension from the previous statement.
\begin{lem}
  \label{lem:ac2}
  For a.e. $0< t_1 < t_2$ and any $\phi\in C_c^\infty(\Omega)$ we
  have
  \begin{equation*}
    \int_{\Gmtt}\phi\,\dmu_{t_2}
    -
    \int_{\Gmto} \phi\,\dmu_{t_1}
    =
    \int_{t_1}^{t_2}\int_{\Gmt} \lan\nabla\phi,\Ht\ran\dmut\,dt
    - \int_{t_1}^{t_2}\int_{\Gmt}\phi |\Ht|^2 \dmut\,dt.
  \end{equation*}
\end{lem}
\begin{proof}
  Let $\phi\in C_c^\infty(\Omega)$ and pick a function $\zeta:\IR\to\IR$
  with compact support and  $\int_0^1\zeta\, dz = 1$. Define
  \begin{equation*}
    \tilde\phi : \Omega\times\IR : (x,z) \to \phi(x)\zeta(z).
  \end{equation*}
  Since $\Ht$ is tangent to $\Omega$, we conclude that
  \begin{equation*}
    \lan \nabla\tilde\phi,\Ht\ran = \zeta \lan\nabla\phi,\Ht\ran.
  \end{equation*}
  Plug $\tilde\phi$ into the statement of lemma \ref{lem:ac1}. 
  As $\tilde\mu = \mu \otimes \CL^1$, is a product and $\tilde\phi$
  is adapted to the product structure, using Fubini's theorem, we can take out
  the integration of $\zeta$ as in the following example
  \begin{equation*}
    \int_\tGmt \tilde\phi\,\tdmu_t
    =
    \int_\IR \int_\Gmt \zeta(z)\phi(x)\, \dmu_t\, dz
    =
    \left(\int_\IR \zeta\,dz \right)
    \,
    \left(\int_\Gmt \phi\,\dmu_t\right)
    =
    \int_\Gmt \phi\,\dmu_t
    .
  \end{equation*}
  This yields the claim.
\end{proof}
We are almost done with the proof of the main theorem \ref{intro.t1},
the only things that remain to be shown is that the statement of
lemma \ref{lem:ac2} holds for all $0<t_1<t_2$ and that we can well
approximate our initial conditions.
\begin{proof}[Proof of Theorem \ref{intro.t1}]
  The general idea will be to approximate arbitrary $t_1$, $t_2$ by
  sequences $t_1^j \geq t_1$ and $t_2^j\geq t_2$ for which by lemma
  \ref{lem:ac1} we have for $\phi \in C_c^\infty(\Omega)$ that
  \begin{equation}
    \label{eq:1}
    \int_{\Gamma_{t_2^j}}\phi\,\dmu_{t_2^j}
    -
    \int_{\Gamma_{t_1^j}}\phi\,\dmu_{t_1^j}
    =
    \int_{t_1^j}^{t_2^j}\int_{\Gmt} \lan\nabla\phi,\Ht\ran\dmut\,dt
    - \int_{t_1^j}^{t_2^j}\int_{\Gmt}\phi |\Ht|^2 \dmut\,dt.
  \end{equation}
  Then we will argue that this statement can be taken to the limit. As
  the function
  \begin{equation*}
    t \mapsto \int_\Gmt \lan \nabla \phi, \Ht \ran - |\Ht|^2 \dmut
  \end{equation*}
  is integrable in $t$, it is clear that the right hand side of
  \eqref{eq:1} converges to
  \begin{equation*}
      \int_{t_1}^{t_2}\int_{\Gmt} \lan\nabla\phi,\Ht\ran\dmut\,dt
    - \int_{t_1}^{t_2}\int_{\Gmt}\phi |\Ht|^2 \dmut\,dt.
  \end{equation*}
  for any sequence $t_1^j\to t_1$ and $t_2^j\to t_2$. The left hand
  side requires a little more argument. To this end, note that as by
  the Brakke flow inequality \ref{intro.3}, we have, as $t_1\leq t_1^j$,
  \begin{equation*}
    |\Gamma_{t_1^j} | \leq |\Gamma_{t_1}|.
  \end{equation*}
  As the characteristic functions $\chi_{K_t}$ are BV functions, and
  since $\chi_{K_{t_1^j}}\to \chi_{K_{t_1}}$ in $L_1$, the lower semi-continuity
  of the total variation of BV-functions implies that
  \begin{equation*}
    |\Gamma_{t_1}| \leq \liminf_{t\to\infty} |\Gamma_{t_1^j}|. 
  \end{equation*}
  Hence 
  \begin{equation*}
    |\Gamma_{t_1}|
    \leq
    \liminf_{t\to\infty} |\Gamma_{t_1^j}|
    \leq
    \limsup_{t\to\infty} |\Gamma_{t_1^j}|
    \leq 
    |\Gamma_{t_1}|,
  \end{equation*}
  and we conclude that $|\Gamma_{t_1^j}|\to |\Gamma_{t_1}|$, as well
  as $|\Gamma_{t_2^j}|\to |\Gamma_{t_2}|$. Now we appeal to lemma
  \ref{lemma:1} and infer that also the left hand side of \eqref{eq:1}
  converges.\\[1.5ex]
Now given an $\Omega\subset \rne$ such that $\p\Omega =: M_0$ is only $C^1$
and carries a nonnegative weak mean curvature in $L^2$ we use that by
lemma \ref{lemma:2} there is a smooth evolution by mean curvature
$M_t, \ 0<t<\gamma$, such that $H_t>0$. We also show in this lemma
that the level set flow of $\p\Omega$ coincides with
the smooth evolution as long as the latter exists. Thus we can do the
whole argument replacing $\Omega$ by $\Omega_t$, where $\Omega_t$ is
the respective open set bounded by $M_t$ for some $t\in
(0,\gamma)$. Using that initially the level set flow is smooth and a
suitable cut-off function we see that (\ref{intro.3}) holds for all
$0\leq t_1 \leq t_2$ and all $\phi \in C^2_c(\rne)$. 
\end{proof}
\begin{lem}
  \label{lemma:1}
  Suppose $E_j\subset \Omega$, $j\geq 1$ and $E\subset \Omega$ are
  Caccioppoli sets, such that $|D\chi_E|(\Omega)<\infty$ and
  $\chi_{E_j}\to \chi_E$ in $L^1(\Omega)$, and
  \begin{equation*}
    \lim_{j\to\infty} |D\chi_{E_j}|(\Omega) = |D\chi_E|(\Omega).
  \end{equation*}
  Then for all $\phi\in C_c(\Omega)$
  \begin{equation*}
    \lim_{j\to\infty} \int_\Omega \phi\, |D\chi_{E_j}|
    =
    \int_\Omega \phi\, |D\chi_E|.
  \end{equation*}  
\end{lem}
\begin{proof}
  We denote $\mu_j = |D\chi_{E_j}|$ and $\mu = |D\chi_E|$.  From
  \cite[Proposition 1.13]{Giusti:1984} we conclude that for every open
  set $A\subset \Omega$ with $\mu(\p A\cap\Omega)=0$, we have
  \begin{equation*}
    \lim_{j\to\infty} \mu_j(A) = \mu(A).    
  \end{equation*}
  Let $A_t := \{\phi > t \}$, then $\p A_t\subset \{ \phi = t \}$,
  whence $\mu(\p A_t \cap \Omega) = 0$ for a.e. $t$. Fix $\eps>0$
  and choose $-T = t_0 < t_1 < \ldots < t_{N_\eps} = T$ such that
  $|t_{i-1}-t_{i}|< \eps$ for $i=1,\ldots,N_\eps$ and $|D\chi_E|(\p
  A_{t_i}\cap\Omega) = 0$ for $i=0,\ldots N_\eps$. Define the step
  function
  \begin{equation*}
    \phi_\eps = t_0 + \sum_{i=1}^{N_\eps} (t_i - t_{i-1})
    \chi_{A_{t_i}}\ .
  \end{equation*}
  It satisfies
  \begin{equation*}
    \sup_\Omega |\phi_\eps - \phi| < \eps    
  \end{equation*}
  and thus
  \begin{equation*}
    \left|\int_\Omega \phi\,\dmu - \int_\Omega \phi_\eps\,\dmu\right|
    \leq \eps \mu(\Omega),
  \end{equation*}
  and
  \begin{equation*}
    \left|\int_\Omega \phi\,\dmu_j - \int_\Omega \phi_\eps\,\dmu_j\right|
    \leq \eps \mu_j(\Omega).
  \end{equation*}  
  Furthermore
  \begin{equation*}
    \begin{split}
      \lim_{j\to\infty} \int_\Omega \phi_\eps\,d\mu_j
      &=
      \lim_{j\to\infty} \Big(t_0 \mu_j(\Omega)
      + \sum_{i=1}^{N_\eps}(t_i-t_{i-1}) \mu_j(A_{t_i})\Big)
      \\
      &= 
      t_0 \mu(\Omega)
      + \sum_{i=1}^{N_\eps}(t_i-t_{i-1}) \mu(A_{t_i})
      =
      \int_\Omega \phi_\eps\,d\mu
    \end{split}
  \end{equation*}
  Thus by letting $\eps\to 0$, we infer the claim.
\end{proof}

In the last lemma we present a slightly stronger version of Lemma 2.6 in
\cite{HuiskenIlmanen02}.

\begin{lem}\label{lemma:2} Let $F_0:M^n\ra \rne$ be a closed, oriented
  hypersurface embedding of class $C^1$, with measurable nonnegative
  weak mean curvature in $L^2(M_0,\hn)$. Then $M_0$ is of class
  $C^1\cap W^{2,2}$ and there exists a smooth evolution by mean
  curvature \mbox{$F:M^n\times(0,\eps)\ra \rne$}$, \ \eps>0,$ such that $M_t \ra M_0$ in
  $C^1\cap W^{2,2}$ and $H_{M_t}>0$ for all $t \in
  (0,\eps)$. Furthermore this smooth evolution coincides with the
  level set flow of $F_0(M)$ as long as the former exists.  
\end{lem}
\begin{proof}
Since $M_0=F_0(M)$ is in $C^1$ and $H \in L^2$ a similar
calculation as in \cite{Giusti:1984}, Appendix A, shows that $M_0$ is
in $W^{2,2}$. Thus by the work of
  Hutchinson, \cite{Hutchinson}, $M_0$ carries a weak second fundamental
  form in $L^2$.
Note that since $M_0$ is compact and in $C^1$ for every $\eps>0$,
there exists a $\delta>0$ such that for all $p\in M_0$,
$$ M_0\cap B_\delta(p) = \text{graph}\, u\ ,$$
where $u:\Omega_p\subset \rn\ra \R,\ \Omega_p$ open, such that
$$ \sup_{\Omega_p}|Du| \leq \eps,\ \ Du(x)=0\ , $$
where $p=(x,u(x))$. By mollification we can pick a sequence of smooth
hypersurfaces $M^i$ converging locally uniformly to $M_0$ in $C^1\cap
W^{2,2}$. The convergence in $C^1$ implies that we can choose for
given $\eps>0$ the above $\delta$ uniform in $i$. Now consider
standard mean curvature flow starting from the approximating surfaces
$F^i:M^n\ra \rne, F^i_0(M)=M^i$. In view of the local gradient
estimates for mean curvature flow in \cite{EckerHuisken91} the surfaces
$M^i_t=F^i(\cdot,t)(M)$ exist on some fixed time interval
$[0,\gamma)$, independent of $i$, and remain controlled graphs in the above
family of coordinate systems relative to $M_0$. Even more by the local
interior estimates in the paper cited above we have
\beq\label{l2.1} \sup_{M^i_t}|A| \leq \frac{C}{t^{1/2}} \qquad t\in
(0,\gamma),
\eeq
independent of $i$, and also uniform estimates, interior in time, on
all higher derivatives. Sending $i\ra \infty$ we extract a limiting
mean curvature flow $M_t,\ \in(0,\gamma)$ satisfying the same
estimates. Note that the uniform local gradient estimate and
(\ref{l2.1}) imply that $M_t \ra M_0$ in $C^{0,\alpha}$ as $t\ra
0$. By the local interior gradient estimates in \cite{EckerHuisken91},
Theorem 2.1, one checks that since $M_0$ is in $C^1$, the surfaces
$M_t$ are equibounded in $C^1$ and thus by Arcela-Ascoli $M_t\ra M_0$
in $C^1$ as $t\ra 0$.\\
Since $M^i\ra M_0$ in $C^1$ and by the local interior gradient
estimates we can now choose for a given $\tilde{\eps}>0$ a 
smooth vectorfield X of unit length on a
$\eta$-neighborhood $U$ of $M_0$ such that, taking $\gamma$ smaller if
necessary,
\beq\label{l2.2}
M^i_t\subset U\ ,\ \ \langle \nu^i(p,t), X(F^i(p,t))\rangle \geq
1-\tilde{\eps}\ ,\ \ |\langle v_{p,t}^i\, ,X(F^i(p,t))\rangle| \leq
\tilde{\eps}\ ,
\eeq
for all $v^i_{p,t} \in T_pM^i_t$ and for all $(p,t) \in M\times
[0,\gamma)$, $ i\geq i_0$. In a local adapted coordinate system we can compute
\beq\label{l2.3}
\frac{d}{dt}\langle X, \nu\rangle =\ \Delta \langle X, \nu\rangle + |A|^2\langle
X, \nu\rangle -2 h^{ij}\langle DX(e_i), e_j\rangle
 -\langle \Delta X, \nu\rangle - H \langle DX(\nu), \nu\rangle\ .
\eeq
If we assume that $\tilde{\eps}<1/4$ we have that $1/2 \leq \langle
X,\nu\rangle -1/4 \leq 3/2$ and we can define
$$ v:=\langle X,\nu\rangle -\frac{1}{4}\ , \qquad w:= \frac{|A|^2}{v^2}\ . $$
Writing the evolution equation of $w$ again in a local adapted coordinate
system, we can estimate
\beq \label{l2.4}\begin{split}
  \frac{d}{dt} w =&\ \Delta w + \frac{4}{v^3}\langle \nabla |A|^2,\nabla v\rangle
  - 6 \frac{|A|^2}{v^4}|\nabla v|^2 - \frac{2}{v^2}|\nabla A|^2
  \\
  &+2\frac{|A|^2}{v^3}\Big( -\frac{1}{4}|A|^2 + 2 h^{ij}\langle DX(e_i), e_j\rangle
  +\langle \Delta X, \nu\rangle + H \langle DX(\nu), \nu\rangle\Big)\,
  \\
  \leq&\ \Delta w + \frac{8|A|}{v^3}|\nabla |A|||\nabla v| -
  \frac{6|A|^2}{v^4}|\nabla v|^2 - \frac{2}{v^2}|\nabla|A||^2
  \\
  &+\frac{|A|^2}{v^3}\Big(-\frac{1}{4}|A|^2 + C(1+|A|)\Big)
  \\
  \leq &\ \Delta w + \frac{2|A|^2}{v^4}|\nabla v|^2 +
  \frac{|A|^2}{v^3}\Big(-\frac{1}{8}|A|^2 + C\Big)\ .
\end{split}\eeq
By (\ref{l2.2}) we can estimate
$$ |\nabla_i v| = |\langle \nabla_i X, \nu\rangle + h_i^{\ j}\langle X,
e_j\rangle| \leq C + \tilde{\eps}|A|\ ,$$
which yields for $\tilde{\eps}$ small enough
$$ \frac{d}{dt} w \leq \Delta w + C w \qquad \text{and} \qquad
\frac{d}{dt}\int_{M^i_t}w \, d\mu \leq C \int_{M^i_t}w\, d\mu\ .$$
Integrating this on $[0,t]$ for $t\leq \gamma$ we see that
$$ \int_{M^i_t}\frac{|A|^2}{v^2}\,
d\mu \leq \exp(Ct) \int_{M^i_0}\frac{|A|^2}{v^2}\,
d\mu \ .$$
Since $M^i_0\ra M_0$ in $W^{2,2}$ this estimate also holds in the
limit. By this estimate $A_t \rightharpoonup A_0$ in $W^{2,2}$, and
since $M_t\ra M_0$ in $C^1$ we have that
$$\lim_{t\ra 0}\int_{M_t}|A|^2\, d\mu = \int_{M_0} |A|^2\, d\mu \ ,$$
which implies full convergence: $M_t \ra M_0$ in $W^{2,2}$. Thus
$(H_t)_- = \min\{H_t,0\} \ra (H_0)_-$ strongly in $L^2$. We can then
check that similar to the computation before the quantity $f:=H/v$
satisfies the evolution equation
$$ \frac{d}{dt}f=\Delta f +\frac{2}{v}\langle \nabla v,\nabla f\rangle
+ \frac{f}{v}\Big(\!\! -\frac{1}{4}|A|^2 + 2 h^{ij}\langle DX(e_i), e_j\rangle
  +\langle \Delta X, \nu\rangle + H \langle DX(\nu), \nu\rangle\Big) $$
and deduce as above that
$$\frac{d}{dt}\int_{M_t}\frac{|H_-|^2}{v^2}\, d\mu \leq C
\int_{M_t}\frac{|H_-|^2}{v^2}\, d\mu \ ,$$
which implies by Gronwall's lemma for $t \in (0,\gamma)$ that
$$ \int_{M_t}\frac{|H_-|^2}{v^2}\, d\mu \leq  \exp(Ct)\int_{M_0}
 \frac{|H_-|^2}{v^2}\, d\mu = 0\ ,$$
proving that $H_t\geq 0$ for $0<t<\gamma$. By the strong maximum
principle and the compactness of $M_t\subset \rne$ it follows that
$H_t>0$ for all $0<t<\gamma$ as required.\\[1ex]
To see that this smooth evolution coincides with the level-set flow of
$M_0$ we define a good coordinate system in a neighborhood of $M_0$. 
Take again $\tilde{M}$ to be
a smooth approximating hypersurface of $M_0$ which is still transverse
to the vectorfield $X$. Let $\Phi_s$ be the flow generated by $X$. Now
define coordinates $\Psi:U\ra \tilde{M}\times (-\eta,\eta)$, where
$U=\bigcup_{-\eta<s<\eta}\Phi_s(\tilde{M})$ for $\eta>0$ small enough,
such that $U$ is a neighborhood of $\tilde{M}$ as follows. We  employ the
flow $\Phi_s$ to 'project' any point $p \in U$ onto
$\tilde{M}_0$ to define the first $n$ coordinates, and the 
parameter $s$ to define the $(n+1)$st coordinate. 
We can  assume that $M_0\subset U$ and thus write $M_0$ in these coordinates as a 'graph' over
$\tilde{M}_0$. Now let $M^s:= \Phi_s(M_0)$ be the translates in
'$x_{n+1}$-direction' in these coordinates and $M^s_t$ be mean curvature
flow with initial condition $M^s$. Since we have locally uniform
gradient bounds in $s$, we can assume that these flows all exist on a
common time interval, say $[0,\varepsilon/2)$, and remain in $U$ for $|s|$ small
enough. Let $u^s:\tilde{M}\times [0, \varepsilon/2) \ra (-\eta,\eta)$ be such
that $\Psi^{-1}(M^s_t)=\text{graph}\big(u^s(\cdot,t)\big)$. Note that
by the interior estimates for mean curvature flow cited before, we
have $u\in C^1(\tilde{M}\times[0,\varepsilon/2))\cap
C^\infty(\tilde{M}\times(0,\varepsilon/2))$. Take $\tilde{g}:=\Psi^*g$ to be the
induced metric on $\tilde{M}\times (-\eta,\eta)$. It can then be checked that the
functions $u^s$ satisfy a parabolic PDE on $\tilde{M}\times
(0,\varepsilon/2)$ of the form
$$D_tu^s =
\bar{g}^{ij}D_{ij}u^s +
f(x,u^s,Du^s)\ , $$
where $\bar{g}^{ij}$ is the inverse of the metric induced on
graph$(u^s)$ by $\tilde{g}$, and $f$ depends smoothly on
$x,u^s,Du^s$. Note that $\bar{g}^{ij}$ depends smoothly on $x,u^s,Du^s$ but
not on $D^2u^s$. By $(\ref{l2.1})$ and the interior gradient estimates
we have that 
$$|Du^s|\leq C\ ,\ \ \ |D^2u^s| \leq \frac{C}{\sqrt{t}}\ ,$$
independent of $s$ for some $C>0$ and all $t\in (0,\varepsilon/2)$. 
Thus interpolating between two solutions $u^{s_1},u^{s_2}$ and
applying the maximum principle we obtain that
$$ \sup_{p\in\tilde{M}} |u^{s_1}(p,t)-u^{s_2}(p,t)| \leq \exp(C
  \sqrt{t}) 
\sup_{p\in\tilde{M}} |u^{s_1}(p,0)-u^{s_2}(p,0)|\ ,
$$
for some constant $C>0$ and all $t\in [0,\varepsilon/2)$. But note
that since the level set flow has to avoid all smooth flows which are
initially disjoint, this implies that for $t\in [0,\varepsilon/2)$ the
level-set flow of $M_0$ coincides with the smooth evolution
$M_t$. Since for $t>0$ the surfaces $M_t$ are smooth it is well-known
that the level-set flow coincides with the smooth evolution as long as the
latter exists.\end{proof} 


\section{The variational principle}\label{varprinc}
As stated in the introduction we give in this final section a variational
formulation for mean curvature flow of mean convex surfaces. Let us
define for $K\subset \Omega$, $K$ compact:

\beq\label{app.1}
J_u(v)=J^K_u(v):=\inl_K |Dv|-\frac{v}{|Du|}\, dx
\eeq

\bd\label{app.d1}
Let $u \in C_\text{loc}^{0,1}(\Omega)\cap L^\infty(\Omega)$ and $|Du|^{-1}\in
L^1_\text{loc}(\Omega)$. Then $u$ is a weak sub- resp. super- solution
of ($\star$) in $\Omega$, if
\beq \label{app.2}
J_u^K(u)\leq J_u^K(v)
\eeq
for every function $v\leq u$ resp. $v\geq u$ which is locally
Lipschitz continuous and satisfies \mbox{$\{v\not = u\}\subset K \subset \Omega$}, where $K$ is
compact.\\[+1ex]
Let $u:\bar{\Omega}\ra [0,\infty)$, $u\in C^{0,1}(\bar{\Omega})$
such that $\{x\in\bar{\Omega}\ |\ u(x)=0\} = \p\Omega$ and
$|Du|^{-1}\in L^1(\Omega)$. Then we call
$u$ a weak solution to $(\star)$ if 
$$ J_u^K(u)\leq J_u^K(v) $$
for every locally Lipschitz continuous function $v$ with 
\mbox{$\{v\not = u\}\subset K \subset \Omega$}, where $K$ is
compact. 

\ed

Since
$$J_u(\min(v,w))+J_u(\max(v,w))=J_u(v)+J_u(w)$$
for $\{v\not = w\}\Subset \Omega$, it follows that $u$ is a weak solution
iff $u$ is a weak sub- and supersolution, provided the boundary conditions
are fulfilled. Note also that the requirement $|Du|^{-1} \in
L^1_\text{loc}$ implies that $u$ is non-fattening,
i.e. $\hne(\{u=t\})=0$ for all $t$.\\[+1.0ex]
\textbf{Equivalent formulation:}
Let $K\subset \Omega$ be compact, and $F\subset \Omega$ be a
Caccioppoli set in a neighborhood of $K$. For a Lipschitz continuous
function $u$ on $\Omega$ with $|Du|^{-1}\in L^1_\text{loc}(\Omega)$ we
can define the functional
\beq\label{app.3}
J_u^K(F):=|\p^*F\cap K|-\inl_{F\cap K}|Du|^{-1}\, dx\ .
\eeq
We say that $E$ minimizes $J_u$
in a set $A$ (from the outside resp. from the inside), if
$$ J_u^K(E)\leq J^K_u(F) $$
for all $F$ with $F\tr E \Subset A$ (with resp. $F\supset E,\ F\subset
E$), with a compact set $K$ with $ F\tr E\subset K \subset
A$.\\[+0.5ex]
By the general inequality
\beq\label{app.4}
J_u(E\cup F)+J_u(E\cap F)\leq J_u(E)+J_u(F)
\eeq
for $E\tr F\Subset A$, it is clear that $E$ minimizes $J_u$ in $A$, 
if it minimizes it from the inside and from the outside. As in
\cite{HuiskenIlmanen} we can show:

\bl\label{app.l1} Let $u\in C^{0,1}_\text{loc}(\Omega)\cap L^\infty(\Omega)$ and
$|Du|^{-1}\in L^1_\text{loc}(\Omega)$. Then $u$ is a weak sub- resp.
supersolution of  $(\star)$ in  $\Omega$, iff for every $t$ the sets $E_t:=\{u>t\}$ minimize
$J_u$ in $\Omega$ from the inside resp. from the outside.
\el

\bpf
1) Let $v$ be  locally Lipschitz continuous with $\{v \not = u\}\subset K\subset \Omega$. For 
\mbox{$F_t:=\{v>t\}$} we have  $F_t\tr E_t \subset K$ for all $t$. Then choose $a<b$ with
$a<u,\!v<b$ on $K$. Using the co-area formula one sees that
\beq\label{app.5}\begin{split}
J_u^K(v)=& \inl_K |Dv|-\frac{v}{|Du|}\, dx\\
=& \int_a^b|\p^*F_t\cap K|\, dt
-\inl_K\int_a^b\chi_{F_t}|Du|^{-1}\, dtdx  - a
\inl_K |Du|^{-1}\, dx\\
=& \int_a^b J_u^K(F_t)\, dt - a\inl_K |Du|^{-1}\, dx\ .
\end{split}
\eeq
Thus, if every $E_t$ minimizes $J_u$ in $\Omega$, then also $u$
minimizes $J_u$. The same works for sub- and supersolutions.\\[+1ex]
2) Now let $u$ be a subsolution of (\ref{app.1}). Choose $t_0$ and
$F$, such that 
$$ F \subset E_{t_0},\ \ E_{t_0}\setminus F \Subset \Omega\ . $$
We aim to show that $J_u(E_{t_0})\leq J_u(F)$. Since $J_u$ is lower semi-continuous
w.r.t. $L^1_\text{loc}$-convergence, we can assume that
\beq\label{app.6} J_u(F)\leq J_u(G) \eeq
for all $G$ with $G \tr E_{t_0} \subset E_{t_0}\setminus F$. Define
\begin{equation*}
F_t:=\begin{cases} F\cap E_t \ \ &t\geq t_0 \\
E_t \ \ \ &0\leq t<t_0
\end{cases}\end{equation*}
By  \eqref{app.6} $J_u(F)\leq J_u(E_t\cup F)$ for all
$t\geq t_0$, and thus by (\ref{app.4})
$$J_u(E_t\cap F)\leq J_u(E_t)$$
for $t\geq t_0$. Thus
$$J_u(F_t)\leq J_u(E_t)\ \ \text{for all}\ t \ .$$
Now define $v$ by $v>t$ on $F_t$, which implies $v\leq u$ and $\{v\not
= u\} \Subset \Omega$. By construction $v \in BV_\text{loc}\cap
L^\infty_\text{loc}$, and $J_u(v)$ is well defined. Approximating $v$
by smooth functions $v_i \rightarrow v$ with $|Dv_i|\rightarrow |Dv|$,
we see that $J_u(u)\leq J_u(v)$ as $u$ is a subsolution. Since then
(\ref{app.5}) also is true for $v$, we have
$$\int_a^bJ_u(E_t)\, dt \leq \int_a^bJ_u(F_t)\, dt\ ,$$
which implies $J_u(F_t)=J_u(E_t)$ for almost all $t$. With
(\ref{app.4}) it follows that
$$J_u(E_t \cup F)\leq J_u(F)$$
for almost all $t\geq t_0$. Taking the limit $t\searrow t_0$ we have
by lower semi-continuity
$$J_u(E_{t_0})\leq J_u(F) \ .$$

3) In the case that $u$ is a supersolution, we choose as in 2) $t_0$ and
$F$ with
$$ E_{t_0}\subset F, \ \ F\setminus E_{t_0}\Subset \Omega\ .$$
As before we can assume that
$$J_u(F)\leq J_u(G)$$
for all $G$ with  $G \tr E_{t_0} \subset F\setminus E_{t_0}$. One
defines again
\begin{equation*}
F_t:=\begin{cases} F\cup E_t \ \ &0\leq t\leq t_0 \\
E_t \ \ \ &t> t_0\ \ ,
\end{cases}\end{equation*}
which leads as above to
$$ J_u(E_t\cap F)\leq J_u(F)$$
for almost all $t\leq t_0$. Since  $|Du|^{- 1}\in
L^1_\text{loc}(\Omega)$, we have \mbox{$\hne(\{u=t_0\})=0$}, and $E_t
\rightarrow E_{t_0}$ for $t\nearrow t_0$. Especially $E_t\cap F \rightarrow
E_{t_0}$, which implies by lower semi-continuity that
$$J_u(E_{t_0})\leq J_u(F)\ .$$
\epf

Applying this equivalent formulation we immediately see that all
superlevelsets of a weak supersolution $u$ minimize area from the outside in $\Omega$:

\bc\label{app.c1}Let $u$ be a weak supersolution on $\Omega$. Then the
sets $E_t$ minimize area from the outside in $\Omega$.
\ec 

To prove uniqueness we aim to show that a weak solution constitutes a
Brakke flow. The idea then is to use the avoidance principle for
Brakke flows to show that the level sets of two weak solutions have to
avoid each other, if they are initially disjoint. In a first step we
show that the mean curvature of almost every level-set is given by
$Du/|Du|^2$, as expected.

\bl\label{app.l1.1} Let $u$ be a weak solution on $\Omega$. Then for
a.e. $t \in [0,T],\ T=\sup_\Omega u$, the sets $\Gamma_t:=\p^*\{u>t\}$,
seen as unit density $n$-rectifiable varifolds, carry a weak mean
curvature $\vec{H}_t \in L^2(\Gamma_t; \hn)$. Furthermore for
a.e. $t$ it holds that
$$ \vec{H}_t=\frac{D u}{|Du|^2}\ .$$
\el 

\bpf Take $X\in C^\infty_c(\Omega;\rne)$ and let $\Phi_s$ be the flow
generated by $X$. We compute:
\beq\label{app.6.1}\begin{split}
0=&\ \frac{d}{ds}\bigg|_{s=0}J_u(u\circ \Phi_s) =
\frac{d}{ds}\bigg|_{s=0}
\int_{-\infty}^{+\infty}\hn(\p^*\{u\circ \Phi_s >  t\})\, dt -
\int_\Omega \frac{u\circ \Phi_s}{|Du|}\, dx \\
=& -\int_{-\infty}^{+\infty}\int_{\Gamma_t} \text{div}_{\Gamma_t}X\,
d\hn\, dt - \int_{\Omega} \frac{\langle Du, X\rangle}{|Du|}\, dx\ .
\end{split}
\eeq
By approximation the last expression still vanishes for any $X \in
C^{0,1}_c(\Omega)$. Let $\Psi:\R \ra \R$ be any smooth function and
replace $X$ above by $\Psi(u)X$. Note that at any point $p\in \Omega$
where $u$ is differentiable and $\Gamma_{u(p)}$ has a weak tangent space
we have $\text{div}_{\Gamma_{u(p)}}(\Psi(u)X) =
\Psi(u)\,\text{div}_{\Gamma_{u(p)}}(X)$. This yields
$$
\int_{-\infty}^{+\infty}\Psi(t)\int_{\Gamma_t}\text{div}_{\Gamma_t}X\,
d\hn\, dt = -\int_{-\infty}^{+\infty}\Psi(t)\int_{\Gamma_t}\Big\langle \frac{Du}{|Du|^2},X\Big\rangle\,
d\hn\, dt\ .
$$
Now let $A$ be a countable dense subset of $C^1_c(\Omega;\, \rne)$. By
the reasoning above there is a set $B\subset [0,T)$ of full measure
such that 
\beq\label{app.6.2} \int_{\Gamma_t}\text{div}_{\Gamma_t}X\, d\hn =
-\int_{\Gamma_t}\Big\langle \frac{Du}{|Du|^2},X\Big \rangle\, d\hn 
\eeq
for all $X\in A$. Since $|Du|^{-1}\in L^1(\Omega)$ we can furthermore
assume that $|Du|^{-1}\in L^2(\Gamma_t,\, \hn)$ and is well-defined
for all $t\in B$. Thus by approximation $(\ref{app.6.2})$ holds for
all $X\in C^1_c(\Omega;\, \rne)$ and $t\in B$. This proves the claim. 
\epf

\bp\label{app.p2} Let $\Omega \subset \rne$ be open and bounded with
$\p\Omega \in C^1$ and $u$ a weak solution of $(\star)$ on
$\Omega$. Then $u$ is a Brakke flow with equality in the sense of
$(\ref{intro.3})$, where $\mu_t:=\hn\res\Gamma_t$ for $t\geq 0$ and
the mean curvature $\vec{H}_t$ of $\Gamma_t$ is given as in Lemma
\ref{app.l1.1}.
\ep
\bpf Take $\phi \in C^\infty_c((0,\infty)),\ \varphi \in
C^1_c(\rne)$, and define a variation $v_s: \Omega \ra \R$ of $u$ by
$$  v_s:= u+s\, \phi(u)\, \varphi\ .$$ 
Note that $\phi(u)\varphi$ has compact support in $\Omega$, and thus $v_s$
is an admissible variation of $u$. We obtain:
\beq\label{app.6.3}\begin{split}
0=&\
\frac{d}{ds}\bigg|_{s=0}J_u(v_s)=\frac{d}{ds}\bigg|_{s=0}\int_{\Omega}
|Dv_s|-\frac{v_s}{|Du|}\, dx\\
=&\ \int_{\Omega}\bigg(\phi^\prime(u)\varphi + \phi(u)\Big(\Big\langle
\frac{Du}{|Du|^2}, D\varphi\Big\rangle -
\frac{\varphi}{|Du|^2}\Big)\bigg)|Du|\, dx\\
=&\ \int_0^\infty\phi^\prime(t)\int_{\Gamma_t} \varphi\, d\hn +
\phi(t) \int_{\Gamma_t} \langle D\varphi, \vec{H}_t\rangle
-\varphi |\vec{H}_t|^2\, d\hn\, dt \ .
\end{split}\eeq
By comparison with small shrinking balls we see that $J_u(\{u>t\})=0$
for all $t\in (0,T)$. Since $u$ is non-fattening this implies that
$\hn(\Gamma_t)\ra 0$ as $t\nearrow T$. Using Lemma \ref{lemma:1} and
the fact that the sets $\{u>t\}$
minimize area from the outside in $\Omega$ we
see that the family of radon measures $\mu_t:=\hn\res \Gamma_t$ is
continuous for $t\geq 0$. Now take $t_1,t_2 \in [0,\infty),\ t_1<
t_2$. Letting $\phi$ appropriately increase to the characteristic
function of the interval $[t_1,t_2]$ we see from $(\ref{app.6.3})$
that $u$ is a Brakke flow with equality as in $(\ref{intro.3})$.
\epf 

\bt\label{app.p3}Let $\Omega \subset \rne$ be open and bounded with
$\p\Omega \in C^1$ and $u$ a weak solution of $(\star)$ on
$\Omega$. Then $u$ is unique.
\et
\bpf 
Let $u_1,u_2$ be two weak solutions to $(\star)$ on $\Omega$. Since
$u_1,u_2 >0$ on $\Omega$ and $\{u_1=0\}=\{u_2=0\}$ we have that
$\{u_2>\tau\}\Subset \{u_1>0\}$ for all $\tau>0$. The avoidance
principle for codimension one Brakke flows, \cite{Ilmanen}, Theorem
10.6, then implies that 
$$ \text{dist}(\Gamma_t^{u_1},\Gamma_{t+\tau}^{u_2}) $$
is increasing in $t$ for all $\tau>0$. Note that Ilmanen's proof of
the avoidance principle also works for the time-integrated version of
a Brakke flow. Since $u_1$ and $u_2$ are continuous this implies that $u_2\leq
u_1$. Repeating this argument with $u_1$ and $u_2$ interchanged we
arrive at the reverse inequality, which implies $u_1=u_2$.
\epf

In the next lemma we show that any smooth mean curvature flow is a
weak sub- and supersolution on the set it sweeps out. To show that
the level set flow is a weak solution on $\Omega$, we later apply this
lemma to the approximating flows $N^{\varepsilon}_t$ on $\Omega\times
\R$ and use corollary \ref{cor:measureconv} to pass to limits.

\bl\label{app.l2}
Let $(N_t)_{c\leq t\leq d}$ be a family of smooth hypersurfaces
$\Omega \times \R$ with strictly positive, uniformly bounded mean
curvature, which flow by mean curvature flow. Let $W$ be the set which
is swept out by the flow
$(N_t)_{c\leq t\leq d}$, and on $W$ let the function $u$ be defined by
$u=t$ on $N_t$ with $E_t:=\{u>t\}$. 
Then the sets $E_t$ minimize $J_u$ on $W$ for all $t\in[c,d]$.
\el
\bpf
The outer unit normal, defined by $\nu_u:=-\frac{Du}{|Du|}$, is a
smooth vector field on $W$ with  div$(\nu_u)=H_{N_t}=|Du|^{-1}>0$\
. For a set $F$ with $F\tr E_t \subset K\Subset W$ we obtain by the
divergence theorem, using $\nu_u$ as a calibration:
\begin{equation*} \begin{split}
|&\p^*E_t\cap K|- \inl_{E_t\cap K}|Du|^{-1}\, dx = \int\limits_{\p E_t \cap K} \nu_{\p E_t}\cdot
\nu_u \, d\hne - \inl_{E_t\cap K}|Du|^{-1}\, dx\\
=& \int\limits_{\p E_t \cap \overline{F}} \nu_{\p E_t}\cdot
 \nu_u \, d\hne +\int\limits_{\p E_t \setminus F} \nu_{\p E_t}\cdot
\nu_u \, d\hne- \inl_{E_t\cap K}|Du|^{-1}\, dx\\
 =& \int\limits_{\p^* F \setminus E_t} \nu_{\p^* F}\cdot
 \nu_u \, d\hne - \inl_{F\setminus E_t}|Du|^{-1}\, dx+\int\limits_{\p^* F \cap \overline{E_t}} 
 \nu_{\p^* F}\cdot \nu_u \, d\hne\\
  &+ \inl_{E_t\setminus F}|Du|^{-1}\, dx- \inl_{E_t\cap K}|Du|^{-1}\, dx\\
=& \int\limits_{\p^* F \cap K} \nu_{\p^* F}\cdot
 \nu_u \, d\hne - \inl_{F\cap K}|Du|^{-1}\, dx\leq |\p^*F\cap K|- \inl_{F\cap K}|Du|^{-1}\, dx\ .
\end{split}\end{equation*}
\epf

\bt\label{app.p1}
Let $\Omega \subset \rne$ be open and bounded. Assume further that
$\p\Omega \in C^1$, carrying a nonnegative weak mean curvature in $L^2$.
Then the level set flow $u:\bar{\Omega} \ra \R$ of $\p\Omega$ is a weak
solution of $(\star)$ on $\Omega$.
\et
\bpf We show that $U((x,z)):=u(x)$, defined on $\Omega\times \R$ is a
weak sub- an supersolution of $(\star)$ on $\Omega\times \R$. That $u$ then is also
a weak sub- and supersolution solution on $\Omega$ follows by a simple
cut-off argument. Note that by Lemma \ref{lemma:2}, $u>0$ on $\Omega$
and so $\{u=0\}=\p\Omega$. \\[+1ex]
We first show that $U$ is a weak supersolution on $\Omega\times \R$. So
take $V\geq U$, $\{U\not = V\} \Subset \Omega\times \R,\ V\in
C^{0,1}_\text{loc}(\Omega\times \R)$. Let $K\subset \Omega\times \R,\ K$ compact with
$\{U\not =V\}\subset K$ and
$$\delta_i:=\max_K|U-U^{\eps_i}|\ ,$$
thus $\delta_i \rightarrow 0$ for $i\rightarrow \infty$. Let
\[V_i:=\begin{cases}\  \max\{U^{\eps_i},V-2\delta_i\}\ \ &\text{for}\ x \in K\\
\ U^{\eps_i}\ \ &\text{for }\ x \not \in K\ ,
\end{cases}\]
We have $V_i\in C^{0,1}_\text{loc}(\Omega\times\R), \ V_i\geq U^{\eps_i},\ \{V_i\not =U^{\eps_i}\}
\subset K$. Furthermore $V_i \rightarrow V$ locally uniformly,
 $V_i=U_i$ on $\Omega\setminus\{V>U\}$ and
\beq \label{app.7}
\mathcal{H}^{n+2}(\{DV_i\not = DV\}\cap \{V>U\})\rightarrow 0\ .
\eeq
By Lemma \ref{app.l2} we have
$$J^K_{U^{\eps_i}}(U^{\eps_i})\leq J^K_{U^{\eps_i}}(V_i)\ ,$$
which can be written as
$$\inl_K |DU^{\eps_i}|+(V_i-U^{\eps_i})|DU^{\eps_i}|^{-1}\, dx\leq
\inl_K|DV_i|\, dx\ .$$
Now we have
$$|DU^{\eps_i}|^{-1} \ra |DU|^{-1}$$
in the sense of Radon measures. 
Since $V_i\rightarrow V$ and $U^{\eps_i} \rightarrow U$ locally uniformly, we see that
$$\inl_K(V-U)|DU|^{-1}\, dx =  \lim_{i\ra \infty} 
\inl_K(V_i-U^{\eps_i})|DU^{\eps_i}|^{-1}\, dx\ .$$
The convergence of $|DU^{\eps_i}|\rhp|DU|$ weakly$-*$ in
$L^\infty(\Omega\times \R)$ as well as
(\ref{app.7}) yield together with the uniform Lipschitz-bound of the
$V_i$'s that
$$ \Big|\!\!\inl_K\!\!|DV_i|-|DV|\, dx\Big| \leq
\Big|\!\!\!\!\!\!\!\!\!\!\!\inl_{\ \ \ K\setminus\{V>U\}}
\!\!\!\!\!\!\!\!\!\!\!\!|DU^{\eps_i}|-|DU|\, dx\Big| + C
\mathcal{H}^{n+2}(\{DV_i\not= DV\}\cap
\{V>U\})\rightarrow 0 \ .$$
Putting this together we have
$$\inl_K |DU|+(V-U)|DU|^{-1}\, dx \leq \inl_K|DV|\, dx\ .$$

That $U$ is also a weak subsolution follows analogously.
\epf 


\bibliographystyle{amsplain}
\bibliography{nodrop} 

\end{document}